\documentclass{article}

\usepackage{amsmath,amssymb,amsthm}
\usepackage{hyperref}
\usepackage{pdfsync}
\usepackage[dvips]{graphicx}
\usepackage{appendix}

\def\R{\mathbb{R}}

\usepackage{float} 
\usepackage{indentfirst} 
\usepackage{float} 

\usepackage{appendix}

\theoremstyle{definition}
\theoremstyle{definition}
\theoremstyle{definition}

\renewcommand{\geq}{\geqslant}
\renewcommand{\leq}{\leqslant}

\renewcommand{\d}{\displaystyle}

\newcommand{\MTCP}{{\bf (MTCP)}}
\newcommand{\OCPZ}{{\bf (OCP0)}}

\newcommand{\OCPReps}{$\bf (OCPR)_\gamma$\ }

\title{Coupled attitude and trajectory optimization for a launcher tilting maneuver}

\author{Jiamin Zhu\footnote{Sorbonne Universit\'es, UPMC Univ Paris 06, CNRS UMR 7598, Laboratoire Jacques-Louis Lions, F-75005, Paris, France (\texttt{zhu@ann.jussieu.fr}).}
\and
Emmanuel Tr\'elat\footnote{Sorbonne Universit\'es, UPMC Univ Paris 06, CNRS UMR 7598, Laboratoire Jacques-Louis Lions, Institut Universitaire de France, F-75005, Paris, France (\texttt{emmanuel.trelat@upmc.fr}).}
\and
Max Cerf\footnote{Airbus Defence and Space, Flight Control Unit, 66 route de Verneuil, BP 3002, 78133 Les Mureaux Cedex, France (\texttt{max.cerf@astrium.eads.net}).}
}

\date{}

\begin{document}
\maketitle

\begin{abstract}
We study the minimum time control problem of the launchers. The optimal trajectories of the problem may contain singular arcs, and thus chattering arcs. The motion of the launcher is described by its attitude kinematics and dynamics and also by its trajectory dynamics. Based on the nature of the system, we implement an efficient indirect numerical method, combined with numerical continuation, to compute numerically the optimal solutions of the problem. Numerical experiments show the efficiency and the robustness of the proposed method.
\end{abstract}

\bigskip


\tableofcontents

\section{Introduction}
The usual way to approach the control of the guidance of satellites consists of considering separately the attitude and the trajectory dynamics. However, for the launchers like Ariane or Pegasus, the trajectories are controlled by their attitude angles, and it is desirable to determine the optimal control subject to the coupled dynamical system. 
Thus, we consider the time minimum control of the attitude reorientation of a rocket taking in consideration both the attitude movement and the trajectory dynamics. We denote this problem as $\MTCP$. 
Our objective is to design an efficient numerical strategy to solve the problem $\MTCP$.
    
There are two main types of numerical approaches to solve an optimal control problem. 
On the one part, the direct methods (see, e.g., \cite{Betts}) consist of discretizing the state and the control and thus of reducing the problem to a nonlinear optimization problem (nonlinear programming) with constraints. 
On the other part, the indirect methods consist of numerically solving a boundary value problem obtained by applying Pontryagin maximum principle (PMP, see \cite{Pontryagin}), by means of a shooting method.
Both direct and indirect methods are not easy to initialize successfully, it is required to combine them with other theoretical or numerical approaches (see the survey \cite{Trelat2}). The numerical continuation is a powerful tool to be combined with the indirect shooting method based on the Pontryagin Maximum Principle. The idea of \emph{continuation} is to solve a problem step by step from a simpler one by parameter deformation (see e.g. \cite{Allgower}). For example, in \cite{CHT,Gergaud,Martinon}, the continuation method is used to solve difficult orbit transfer problems. Here, we will use numerical continuation to numerically solve the problem.    
    
Indeed, the coupling of the attitude movement (fast) with the trajectory dynamics (slow) in the problem $\MTCP$ generates difficulties in numerical approaches, especially in the indirect methods, where the Newton-like methods are used to solve the boundary value problem. Moreover, this property makes the numerical approaches extremely difficult to be initialized.
Another difficulty in numerically solving the problem $\MTCP$ is the chattering phenomenon. \emph{Chattering} means that the control switches an infinite number of times over a compact time interval. Such a phenomenon typically occurs when trying to connect bang arcs with a higher-order singular arc (see, e.g., \cite{FULLER1, Marchal, ZELIKIN}). We proved in \cite{ZTC1,ZTC2} that there exists such a chattering phenomenon in the problem $\MTCP$, and the chattering phenomenon is a bad news when applying numerical approaches.

In view of these difficulties, we propose an efficient numerical continuation procedure to compute optimal and sub-optimal trajectories for the problem $\MTCP$. Due to the chattering phenomenon, numerical continuation combined with shooting cannot give an optimal solution to the problem for certain terminal conditions for which the optimal trajectory contains a singular arc of higher-order. In that case, our indirect approach generates sub-optimal solutions, by stopping the continuation procedure before its failure due to chattering. 

Our objective is to design and implement a numerical method having the following qualities: general (vehicle features, terminal conditions), robust, fast and automatic, so that it could be used efficiently by an inexperienced user on a large scope of applications of the real world.
Numerical experiments indicates that this approach happens to meet to all the expected features.

\medskip

The paper is organized as follows. In Section \ref{pbstate}, we establish the model of the rocket movement and we formulate the problem $\MTCP$. In Section \ref{Theres}, the PMP is applied to the problem $\MTCP$ and some theoretical results on the extremals are introduced. In Section \ref{AuxiPro}, two auxiliary problems are introduced and the continuation procedure is stated. Some numerical tips improving the proposed numerical strategy are also presented. Finally in Section \ref{Chp_numerical}, numerical results are given.

\section{Problem Statement} \label{pbstate}
We will consider two types of launchers: one is a classical launcher like Ariane rocket, the other one is an airborne launcher like Pegasus rocket (being launched horizontally from an airplane). In the case of an Ariane-type launcher, the aerodynamical forces are very small compared with its thrust and the gravity, and so we ignore them. However, in case of a Pegasus-type launcher, we will need to add the aerodynamical forces to the model. 
Throughout the paper, we make the following assumptions:
\begin{itemize}
\item The Earth is a sphere and is fixed in the inertial space.
\item The position and the mass of the rocket are constant during the maneuver.
\item The rocket is an axial symmetric cylinder.
\item The rocket engine cannot be shut off during the maneuver and the module of the thrust force is constant, taking its maximum value, i.e., $T=T_{max}$.
\end{itemize}

\subsection{Coordinates and Model}
All the coordinate systems introduced here are Cartesian coordinate systems.

The \textbf{launch frame} (reference frame) $S_{R}=(\hat{x}_R,\hat{y}_R,\hat{z}_R)$ is fixed around the launch point $O_R$ (where the rocket is launched). The axis $\hat{x}_R$ is pointing radially outwardly (normal to the local tangent plane), and the axis $\hat{z}_R$ points to the North.

The \textbf{body frame} $S_b=(\hat{x}_b,\hat{y}_b,\hat{z}_b)$ is defined as follows. The origin of the frame $O_b$ is fixed around the mass center of the rocket, the axis $\hat{z}_b$ is along the axis-symmetric axis of the rocket, and the axis $\hat{x}_b$ is in the cross-section. The body frame can be derived by three ordered unit single-axis rotations from the launch frame, 
$$
    S_{R} \xrightarrow{R_y(\theta)} \circ \xrightarrow{R_x(\psi)} \circ \xrightarrow{R_z(\phi)} S_b,
$$
where $\theta$ is the pitch angle, $\psi$ is the yaw angle, $\phi$ is the roll angle and $R_a(b)$ means to rotate the frame around the axis $a \in \{x,y,z\}$ by an angle $b \in \mathbb{R}$. Therefore, the transfer matrix from $S_R$ to $S_b$ is $L_{bR}= R_z(\phi) R_x(\psi) R_y(\theta)$. 

The \textbf{velocity frame} $S_v=(\hat{x}_v,\hat{y}_v,\hat{z}_v)$ is fixed around the mass center of rocket. The axis $\hat{x}_v$ is parallel to the velocity of the rocket $\vec{v}$, and the axis $\hat{z}_v$ is normal to the velocity, pointing to the direction of the lift force $\vec{L}$ of the rocket. This frame can be derived by two unit single-axis rotations from the launch frame, 
$$
    S_{R} \xrightarrow{R_x(\xi)} \circ \xrightarrow{R_y(\kappa)}  S_v,
$$
where $\xi$ is the flight path angle and $\kappa$ is the bank angle of the flight. Denote by $v$ the module of the velocity vector, i.e., $v=\sqrt{v_x^2+v_y^2+v_z^2}$.
since 
$(\vec{v})_R
= (v \cos \xi, v \sin \xi \sin \kappa, -v \sin \xi \cos \kappa)^\top$,
we have 
$$
\cos \xi = v_x/v, \quad \tan \kappa = -v_y/v_z.
$$
This frame will be used to introduce aerodynamical forces (lift $\vec{L}$ and drag $\vec{D}$) into the model.

\paragraph{Attitude motion}
The dynamical and kinematic equations of the attitude are
\begin{equation}\label{sys_atti}
\begin{split}
& \dot{\theta}=(\omega_x \sin \phi + \omega_y \cos \phi)/ \cos \psi , \quad \dot{\psi}=\omega_x \cos \phi - \omega_y \sin \phi ,\quad \dot{\phi}= (\omega_x \sin \phi + \omega_y \cos \phi) \tan \psi ,\\
& \dot{\omega }_x= -\bar{b} u_2,\qquad \dot{\omega }_y=  \bar{b} u_1.
\end{split}
\end{equation}
where $u=(u_1,u_2) \in \mathbb{R}^2$ is the control input expressed in frame $S_b$, $(\vec{\omega})_b = (\omega_x,\omega_y)$ is the angular velocity of the launcher $\vec{\omega}$ expressed in frame $S_b$ (assuming that there is no rotation along the asymptotical axis of the launcher), and $\bar{b}$ is the angular acceleration constant depending on $\mu_{max}$ (the maximum angle between the thrust $\vec{T}$ and the axis $\hat{z}_b$), $I_c$ (the inertia around the axis $\hat{x}_b$), $l$ (the length of the rocket). Note that $\| u \| = \sqrt{u_1^2+u_2^2} \leq 1$. See \cite{ZTC2} for details of this model.

\paragraph{Trajectory dynamics}
The drag $\vec{D}$ and lift $\vec{L}$ are calculated by $
\vec{D} = - \frac{1}{2} \rho v^2 S C_x \hat{x}_v$, $\vec{L} =  \frac{1}{2} \rho v^2 S C_x \hat{z}_v$, where $\rho$ is the air density, $S$ is the reference surface, $C_x$ and $C_z$ are constant coefficients of the drag and the lift. In the velocity frame $S_v$, the components of the vectors $\vec{D}$ and $\vec{L}$ are expressed as $(\vec{D})_v = (-c_x m v^2,0,0)^\top$, $(\vec{L})_v = (0,0,c_z m v^2)^\top$ where $c_x = \frac{1}{2} \rho S C_x/m$ and $c_z=\frac{1}{2} \rho S C_x/m$. By using the transfer matrix from $S_v$ to $S_R$, we get 
$$
(\vec{D})_R = c_x m v^2(- \cos \xi,- \sin \xi \sin \kappa,  \sin \xi \cos \kappa)^\top, \quad 
(\vec{L})_R = c_z m v^2( - \sin \xi, - \cos \xi \sin \kappa, \cos \xi \cos \kappa)^\top,
$$
and so we get the trajectory dynamical equations in frame $S_R$ as
\begin{equation}\label{sys_orbit_pegasus}
\begin{split}
& \dot{v}_x= a \sin \theta \cos \psi + g_x - c_x v^2\cos \xi-c_z v^2\sin \xi, \\
& \dot{v}_y= - a \sin \psi + g_y - c_x v^2\sin \xi \sin \kappa-c_z v^2 \cos \xi \sin \kappa,\\
& \dot{v}_z= a \cos \theta \cos \psi + g_z + c_x v^2\sin \xi \cos \kappa+c_z v^2\cos \xi \cos \kappa,
\end{split}
\end{equation}
where $(\vec{v})_R = (v_x,v_y,v_z)$ is the velocity of the launcher $\vec{v}$ expressed in $S_R$, $(\vec{g})_R=(g_x,g_y,g_z)$ is the gravity vector expressed in $S_R$, $a = T_{max}/m$ is the thrust acceleration. For the Ariane rocket, we do not consider the aerodynamical forces, i.e., $c_x=c_z=0$.

\subsection{Minimum Time Control Problem $\MTCP$}
Defining the state variable $x$ and the control variable $u$ as
$$
x=(v_x, v_y, v_z, \theta, \psi, \phi, \omega_x, \omega_y), \quad u=(u_1,u_2).
$$
The initial conditions are defined by ${v_{x_0}}$, ${v_{y_0}}$, ${v_{z_0}}$, $\theta_0$, $\psi_0$, $\phi_0$, ${\omega_{x_0}}$ and ${\omega_{y_0}}$,
\begin{equation} \label{OCPc_ic}
\begin{split}
& v_x(0) = {v_{x_0}},\quad  
v_y(0) = {v_{y_0}},\quad  
v_z(0) = {v_{z_0}},\\  
& \theta(0) = \theta_0,\quad \psi(0)=\psi_0,\quad \phi(0)=\phi_0,\quad
\omega_x(0) = {\omega_{x_0}},\quad \omega_y(0)={\omega_{y_0}}.
\end{split}
\end{equation}

The desired final velocity is required to be parallel to the axis $\hat{z}_b$, according to $
(\vec{V}(t_f))_R \wedge (\hat{z}_b (t_f))_R=\vec{0}$. The constraints on the final conditions are defined by $\theta_f$, $\psi_f$, $\phi_f$, $\omega_{x_f}$ and $\omega_{y_f}$.
\begin{equation} \label{OCPc_fc}
\begin{split}
& v_{z_f} \sin \psi_f + v_{y_f} \cos \theta_f \cos \psi_f =0,\quad v_{z_f} \sin \theta_f - v_{x_f} \cos \theta_f =0,\\
& \theta(t_f)=\theta_f,\quad \psi(t_f)=\psi_f,\quad \phi(t_f)=\phi_f, \quad \omega_x(t_f)=\omega_{x_f},\quad \omega_y(t_f)=\omega_{y_f}.
\end{split}
\end{equation}
Note that the parallel condition on the final velocity is due to the fact that most rockets are planned to maintain a zero angle of attack along the flight. The angle of attack, when the wind is set to zero, is defined as the angle between the velocity and the rocket body axis.

We set $x_0 = ({v_{x_0}},{v_{y_0}},{v_{z_0}},\theta_0,\psi_0,\phi_0,{\omega_{x_0}},{\omega_{y_0}})$, and we define the target set
\begin{equation*}
\begin{split}
M_1 = & \{(v_x,v_y,v_z,\theta,\psi,\phi,\omega_x,\omega_y) \in \R^8 \ \mid\  v_z \sin \psi_f + v_y \cos \theta_f \cos \psi_f =0,\\
 &\qquad v_z \sin \psi_f + v_y \cos \theta_f \cos \psi_f =0, \quad \theta =\theta_f,\quad \psi=\psi_f,\quad \phi=\phi_f,\quad \omega_x=\omega_{x_f},\quad \omega_y=\omega_{y_f} \} .
\end{split}
\end{equation*}
The minimum time control problem $\MTCP$ consists of steering the system \eqref{sys_atti}-\eqref{sys_orbit_pegasus} from $x(0)=x_0$ to the final target $M_1$ in minimum time $t_f$, with controls satisfying the constraint $u_1^2+u_2^2 \leq 1$.

\medskip

For convenience, the velocity will be defined by polar coordinates, i.e., $v_x = v \sin \theta_v \cos \psi_v$, $v_y = -v \sin \psi_v$, $v_z = v \cos \theta_v \cos \psi_v$, where $\theta_v$ and $\psi_v$ are ``pitch" and ``yaw" angles of the velocity vector.

\section{Theoretical Results} \label{Theres}

According to the Pontryagin Maximum Principle (PMP), the Hamiltonian of the optimal control problem $\MTCP$ is 
$$
H(x(t),p(t),p^0,u(t)) = \langle p(t), f(x(t),u(t)) \rangle + p^0,
$$
where $p=(p_{vx},p_{vy},p_{vz},p_{\theta},p_{\psi},p_{\phi},p_{\omega x},p_{\omega y})$ is the adjoint variable, $f(x,u)$ is the derivative of the state variable $x$, i.e., $\dot{x}=f(x,u)$ and $p^0 \leq 0$ is a real number. Here we assume $p^0=-1$ (see \cite{ZTC2}).
The differential equation of the adjoint variable $p$ is given by the partial derivative of the Hamiltonian as
\begin{equation}
\label{adjointsys}
\dot{p}(t) = -\frac{\partial H}{\partial x}(x(t),p(t),p^0,u(t)).
\end{equation}

The maximization condition of the PMP yields, almost everywhere on $[0,t_f]$,
\begin{equation*} 
u(t) = \frac{(h_1(t),h_2(t)) }{\sqrt{h_1(t)^2+h_2(t)^2}} = \frac{\Phi(t)}{\Vert\Phi(t)\Vert},
\end{equation*}
whenever $\Phi(t)=(h_1(t),h_2(t))\neq (0,0)$. Here we have $h_1(t)  =\bar{b} p_{\omega y}(t)$ and $h_2(t) = = - \bar{b} p_{\omega x}(t)$.
We call $\Phi$ (as well as its components) the switching function.

Moreover, we have the transversality condition $p(t_f) \perp T_{x(t_f)} M_1$, where $T_{x(t_f)}M_1$ is the tangent space to $M_1$ at the point $x(t_f)$, i.e.,
\begin{equation}
\label{condtransv}
p_{v_y} \sin \psi_f = p_{v_x} \sin \theta_f \cos \psi_f+ p_{v_z} \cos \theta_f \cos \psi_f,
\end{equation}
and, the final time $t_f$ being free and the system being autonomous, we have also 
$h_0(x(t),p(t))+\Vert\Phi(t)\Vert+p^0=0,\: \forall t\in[0,t_f]$.

The quadruple $(x(\cdot),p(\cdot),p^0,u(\cdot))$ is called an extremal lift of $x(\cdot)$.
An extremal is said to be normal (resp., abnormal) if $p^0 < 0$ (resp., $p^0 = 0$). 
We say that an arc (restriction of an extremal to a subinterval $I$) is \emph{regular} if $\Vert \Phi(t)\Vert \neq 0$ along $I$. Otherwise, the arc is said to be \emph{singular}.

A \emph{switching time} is a time $t$ at which $\Phi(t)=(0,0)$, that is, both $h_1$ and $h_2$ vanish at time $t$.
An arc that is a concatenation of an infinite number of regular arcs over a compact time interval is said to be \emph{chattering}. The chattering arc is associated with a \emph{chattering control} that switches an infinite number of times, over a compact time interval.
A junction between a regular arc and a singular arc is said to be a \emph{singular junction}.

When $c_x=0$ and $c_z=0$, the explicit form of equation \eqref{adjointsys} is given by
\begin{equation} \label{sys_adjoint}
\begin{split}
\dot{p}_{v_x}&= 0,\quad \dot{p}_{v_y}= 0,\quad \dot{p}_{v_z}= 0,\\
\dot{p}_{\theta}&= -a \cos \psi (\,p_{v_x} \cos \theta -\,p_{v_z} \sin \theta),\\
\dot{p}_{\psi}&= a \sin \psi \sin \theta \,p_{v_x} + a \cos \psi \,p_{v_y} + a \cos \theta \sin \psi \,p_{v_z}
-\sin \psi (\omega_x \sin \phi + \omega_y \cos \phi)/\cos^2 \psi \,p_{\theta} \\
&\qquad\qquad\qquad  -(\omega_x \sin \phi + \omega_y \cos \phi)/ \cos^2 \psi \,p_{\phi}, \\
\dot{p}_{\phi}& = -(\omega_x \cos \phi -\omega_y \sin \phi)/\cos \psi \,p_{\theta}
+(\omega_x \sin \phi + \omega_y \cos \phi) \,p_{\psi} \\
&\qquad\qquad\qquad -\tan \psi (\omega_x \cos \phi -\omega_y \sin \phi) \,p_{\phi}, \\
\dot{p}_{\omega_x}&= - \sin \phi / \cos \psi \,p_{\theta}- \cos \phi \,p_{\psi}
 -\sin \psi \sin \phi / \cos \psi \,p_{\phi}, \\
\dot{p}_{\omega_y}&=  - \cos \phi / \cos \psi \,p_{\theta}+ \sin \phi \,p_{\psi}
 -\sin \psi \cos \phi / \cos \psi \,p_{\phi} .
\end{split}
\end{equation}
We have shown in \cite{ZTC2} that the singular extremals of the problem $\MTCP$ are normal ones, i.e., $p^0 \ne 0$, and that they are of intrinsic order two. In that case the singular junction can only be realized by chattering, i.e., if the singular junction happens at time $\tau$ and the control is regular on $(t_1,\tau)$ and is singular on $(\tau,t_2)$, then the control has to switch an infinite number of times on $(\tau-\epsilon,\tau)$, $\forall \epsilon \in \mathbb{R}$ such that $0< \epsilon < t_1-\tau$.

The concatenation of regular arcs are classified in \cite{ZTC2} by their contact with the switching surface (filled by switching points), and we know that if the switching points of the problem $\MTCP$ are of order one, then the control will turn an angle $\pi$ when passing the switching surface. We will see this phenomenon in the numerical results.

Actually, when $c_x \ne 0$ and $c_z \ne 0$, the derivative of $p_v=(p_{vx},p_{vy},p_{vz})$ is no longer zero. However, the terms induced by the air draft and lift do not change the Lie bracket configuration of the system, and thus the results obtained in \cite{ZTC2} are still valid. 

\medskip

Note again that the chattering causes hard problems when applying numerical methods. In that case, we can no longer obtain an optimal solution numerically. However, a sub-optimal solution can be derived by stopping the continuation procedure proposed in Section \ref{indirecte_strategy} before chattering occurs.

\section{Auxiliary Problems and Continuation Procedure} \label{AuxiPro}

To construct a suitable continuation procedure, we introduce two auxiliary problems. 
The first one is the \emph{problem of order zero}, denoted by $\OCPZ$. The solution of this problem can be easily computed. Then we just need to plug this simple, low-dimensional solution in higher dimension, in order to initialize an indirect method for the more complicated problem $\MTCP$.
The second one is the \emph{regularized problem}, which is used to overcome the chattering issues. Since the solution of problem $\OCPZ$ is contained in the singular surface (filled by the singular solutions) for the problem $\MTCP$, passing directly from the solution of problem $\OCPZ$ to the problem $\MTCP$ makes the optimal extremals to contain a singular arc (and thus chattering arcs), and the shooting method will certainly fail due to the numerical integration of discontinuous Hamiltonian system.

\subsection{Two auxiliary problems}
\label{pb_auxi}
\paragraph{Problem of order zero.}
We define the \emph{problem of order zero}, as a ``subproblem'' of the complete problem $\MTCP$, in the sense that we consider only the trajectory dynamics (without aerodynamical forces) and that we assume a perfect control meaning that
 the attitude angles (Euler angles) take the desired values instantaneously. Thus, the attitude angles are considered as control inputs in that simpler problem. Denoting the rocket axial symmetric axis as $\vec{e}$ and considering it as the control vector (which is consistent with the attitude angles $\theta$, $\psi$), we formulate the problem as follows:
$$
\dot{\vec{v}} = a \vec{e}+\vec{g}, \quad
\vec{v}(0)=\vec{v}_0,\quad \vec{v}(t_f) // \vec{w},\quad  \Vert \vec{w} \Vert =1,  \qquad \min t_f ,
$$
where $\vec{w}$ is a given vector that refers to the desired target velocity direction. 
This problem is easy to solve, and the solution derived by applying the PMP is
\begin{equation*}
 \vec{e}^{\ast}=\frac{1}{a}\left(\frac{k \vec{w}-\vec{v}_0}{t_f}-\vec{g}\right), \quad
 t_f=\frac{-a_2 + \sqrt{a_2^2-4a_1a_3}}{2 a_1},\quad
 \vec{p}_v=\frac{-p^0}{a + \langle \vec{e}^{\ast},\vec{g} \rangle } \vec{e}^{\ast}.
\end{equation*}
with 
$k=\langle \vec{v}_0,\vec{w} \rangle + \langle \vec{g},\vec{w} \rangle t_f$, 
$a_1 = a^2-\Vert  \langle \vec{g},\vec{w} \rangle \vec{w}-\vec{g} \Vert ^2$, 
$a_2 = 2 ( \langle \vec{v}_0,\vec{w} \rangle \langle \vec{g},\vec{w} \rangle - \langle \vec{v}_0,\vec{g} \rangle)$, and
$a_3 = - \Vert  \langle \vec{v}_0,\vec{w} \rangle \vec{w}-\vec{v}_0 \Vert ^2$.

Since the vector $\vec{e}$ is expressed in the launch frame $S_r$ as $(\vec{e})_R = (\sin \theta \cos \psi, -\sin \psi, \cos \theta \cos \psi)^{\top} $, 
the Euler angles $\theta^{\ast} \in (-\pi,\pi)$ can be calculated by assuming $\psi^\ast \in (-\pi/2,\pi/2)$ as
\begin{equation} \label{sing_angles_tet}
\theta^{\ast} = \mathrm{atan2}(e_x^\ast,e_z^\ast) = \mathrm{sign}(e_x^\ast)
\begin{cases}
\mathrm{arctan}(| e_x^\ast/e_z^\ast |), \,\, e_z^\ast >0,\\
\pi/2, \,\, e_z^\ast =0,\\
(\pi - \mathrm{arctan}(| e_x^\ast/e_z^\ast |)), \,\, e_z^\ast <0, 
\end{cases}
\end{equation}
where $\mathrm{sign}(\cdot)$ is the sign function, and
\begin{equation} \label{sing_angles_psi}
\psi^{\ast}= \mathrm{arcsin} (-e_y^\ast).
\end{equation}
The assumption of $\psi^\ast \in (-\pi/2,\pi/2)$ generally works for the Ariane-type rockets since the maneuvers are generally quasi planar and the change in the yaw angle is small. 
However, for the Pegasus-type rockets, which could make large yaw angle maneuvers, we have to consider also the case $\cos \psi^\ast<0$ leading to
\begin{equation}\label{sing_angle_2}
\theta^{\ast} = \mathrm{atan2}(- e_x^\ast, - e_z^\ast) ,\quad
\psi^{\ast} = - \mathrm{sign}(e_y^\ast) (\pi -  \mathrm{arcsin} (| -e_y^\ast|)).
\end{equation}
Both \ref{sing_angle_2} and \eqref{sing_angles_tet} and \eqref{sing_angles_psi} are possible solutions. We will choose the pair $(\theta^\ast,\psi^\ast)$ minimizing the value $|\psi_0-\psi^\ast|+|\psi_f-\psi^\ast|$. 

Indeed, when $\psi = \pm \pi/2 + k \pi$, $k \in \mathbb{Z}$, the Euler angles are not well defined, and if $\psi_0$ and $\psi_f$ are given near $\pm \pi/2$, we need to ``change" them to a be near $\pm \pi/2$. This is done in section \ref{precond} below by transforming of the reference frame $S_R$ to a new reference frame $S_R^\ast$. The singularities of Euler angles are also treated in section \ref{SingularEuler} by calculating the limit values of the vector field at these singular points.

\medskip

Given a real number $\phi^\ast$, the optimal solution of the problem $\OCPZ$ actually corresponds to a singular solution of $\MTCP$ with the terminal conditions given by
\begin{equation} \label{inicond_ocpztomtcp}
\begin{split}
v_x(0) = {v_{x_0}}, \quad v_y(0) = {v_{y_0}}, \quad v_z(0) = {v_{z_0}}, \\
\theta(0)=\theta^\ast, \quad \psi(0)=\psi^\ast,\quad, \phi(0)=\phi^\ast,\quad \omega_x(0) = 0,\quad \omega_y(0)=0,
\end{split}
\end{equation}
\begin{equation} \label{fincond_ocpztomtcp}
v_z(t_f)\sin \psi_f + v_y(t_f) \cos \theta_f \cos \psi_f=0,\quad v_z(t_f)\sin \theta_f-v_x(t_f)\cos \theta_f=0,
\end{equation}
\begin{equation} \label{fincond_ocpztomtcp2}
 \theta(t_f)=\theta^\ast, \quad \psi(t_f)=\psi^\ast,\quad, \phi(t_f)=\phi^\ast,\quad \omega_x(t_f) = 0,\quad \omega_y(t_f)=0.
 \end{equation}
It is worth nothing that this solution is lying in the singular surface of the problem $\MTCP$ meaning that it is on the ``highway" between two fixed points in the state space. On this ``highway", the system goes the most rapidly towards to aimed final point or manifold. Indeed, we observe in the numerical simulations that the singular extremals of the problem $\MTCP$ acts a role similar to that of the stable points in the ``turnpike" phenomenon as described in \cite{TrelatZuazua_JDE2015}: the optimal trajectories first tend to reach the singular surface (to have a greater speed in transferring its state), then stay in the singular surface for a while until it is sufficiently close to the target submanifold $M_1$, and finally get off the singular surface to join the target submanifold. Note that the singular arc is not necessary if the state is sufficiently close to the target.

In view of \eqref{inicond_ocpztomtcp}-\eqref{fincond_ocpztomtcp2}, a natural continuation strategy is to simply vary the terminal conditions step by step until they correspond to $M_1$. 
However, as we have mentioned, the chattering phenomenon causes the failure of this simple strategy, and thus we introduce another auxiliary problem, the regularized problem, in which the singular arcs of the problem $\MTCP$ become regular arcs, and thus the difficulty caused by chattering is bypassed.

\paragraph{Regularized problem.}
Let $\gamma > 0$ be arbitrary. The regularized problem \OCPReps consists of minimizing the cost functional
\begin{equation} \label{cost_opcr}
    C_\gamma = t_f + \gamma \int_0^{t_f} (u_1^2 + u_2^2) \, dt ,
\end{equation}
for the system \eqref{sys_atti}-\eqref{sys_orbit_pegasus}, under the control constraints $-1\leq u_i\leq 1$, $i=1,2$, with terminal conditions \eqref{OCPc_ic}-\eqref{OCPc_fc}.
Note that, here, we replace the constraint $u_1^2+u_2^2\leq 1$ with the constraints $-1\leq u_i\leq 1$. The advantage, for this intermediate optimal control problem with the cost \eqref{cost_opcr}, is that the extremal controls are continuous.

The Hamiltonian is
\begin{equation}\label{Hamil_n}
H_\gamma=\langle p,f(x,u) \rangle +p^0(1+\gamma u_1^2 + \gamma u_2^2),
\end{equation}
and according to the PMP, the optimal controls are
\begin{equation} \label{OCPn_u1}
u_1(t) = \mathrm{sat}(-1,-\bar{b} p_{\omega_y}(t) / (2 \gamma p^0),1),\qquad
u_2(t) = \mathrm{sat}(-1,\bar{b} p_{\omega_x}(t) / (2 \gamma p^0),1),
\end{equation}
where the saturation operator $\mathrm{sat}$ is defined by $\mathrm{sat}(-1,f(t),1)=-1$ if $f(t)\leq -1$; $1$ if $f(t)\geq 1$; and $f(t)$ if $-1\leq f(t)\leq 1$.

\medskip
 
Indeed, an extremal of $\OCPZ$ can be embedded into the problem \OCPReps by setting
$$
u(t) = (0,0), \quad \theta(t)=\theta^\ast, \quad \psi(t)=\psi^\ast, \quad \phi(t)=\phi^\ast, \quad \omega_x(t) = 0, \quad \omega_y(t)=0,
$$
$$
p_{\theta}(t)=0, \quad p_{\psi}(t)= 0, \quad p_{\phi}(t)=0, \quad p_{\omega x}(t) = 0, \quad p_{\omega y}(t)=0,
$$
where $\theta^\ast$ and $\psi^\ast$ are given by \eqref{sing_angles_tet} and \eqref{sing_angles_psi}, 
with terminal conditions given by \eqref{inicond_ocpztomtcp}-\eqref{fincond_ocpztomtcp2}. Moreover, it is not a singular extremal for the problem \OCPReps.
The extremal equations for \OCPReps are the same than for $\MTCP$, as well as the transversality conditions.

Note that the difficulty caused by the chattering still exists when trying to go from the regularized problem back to the problem $\MTCP$ if there exists a singular arc in the optimal trajectory. However, in that case, by stopping at some point during the last continuation step (see below), an acceptable sub-optimal trajectory will be found for the problem $\MTCP$.

\subsection{Continuation procedure} \label{indirecte_strategy}
The ultimate objective is to compute the optimal solution of the problem $\MTCP$, starting from the explicit solution of $\OCPZ$. 
\paragraph{Numerical strategy}
We proceed as follows:
\begin{itemize}
\item First, we embed the solution of $\OCPZ$ into \OCPReps.
For convenience, we still denote by $\OCPZ$ the problem $\OCPZ$ seen in high dimension. 
\item Then, we pass from $\OCPZ$ to $\MTCP$ by means of a numerical continuation procedure, involving four continuation parameters: the first three parameters $\lambda_1$, $\lambda_2$ and/or $\lambda_3$ are used to pass continuously from the optimal solution of $\OCPZ$ to the optimal solution of the regularized problem \OCPReps, for some fixed $\gamma>0$, and the last parameter $\lambda_4$ is then used to pass to the optimal solution of $\MTCP$ (see Figure \ref{homotopie}).
\end{itemize}

\begin{figure}[h]
\centering 
	\includegraphics[scale=0.8]{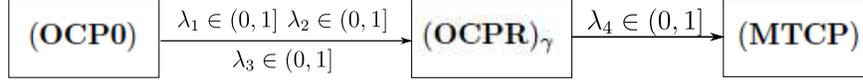}
	\caption{Continuation procedure.}
	\label{homotopie}
\end{figure}

The unknowns of the shooting problem are $p_{v_x}(0)$, $p_{v_y}(0)$, $p_{v_z}(0)$, $p_{\theta}(0)$, $p_{\psi}(0)$, $p_{\phi}(0)$, $p_{\omega_x}(0)$, $p_{\omega_y}(0)$ and $t_f$. When $D=0$ and $L=0$, we have that $p_{v_x}$, $p_{v_y}$ and $p_{v_z}$ are constants, hence by using the transversality condition \eqref{condtransv},
we can easily reduce the number of the unknowns by one. In view of more general cases, we will not do this reduction and we will use the transversality condition as a part of the shooting function. In the following, we denote the continuation step corresponding to parameter $\lambda_i$ by $\lambda_i$-continuation.

\paragraph{$\lambda_1$-continuation and $\lambda_2$-continuation}
The parameter $\lambda_1$ is used to act, by continuation, on the initial conditions, according to
$$
\theta(0) = \theta^{\ast} (1-\lambda_1) + \theta_0 \lambda_1, \quad \psi(0) = \psi^{\ast} (1-\lambda_1) + \psi_0 \lambda_1, \quad \phi(0) = \phi^{\ast} (1-\lambda_1) + \phi_0 \lambda_1, 
$$
$$
\omega_x(0) = \omega_x^{\ast} (1-\lambda_1) + {\omega_{x_0}} \lambda_1, \quad \omega_y(0) = \omega_y^{\ast} (1-\lambda_1) + {\omega_{y_0}} \lambda_1,
$$
where $\omega_x^{\ast} = \omega_y^{\ast} =0$, $\phi^{\ast} = 0$, and $\theta^{\ast}$, $\psi^{\ast}$ are calculated through equation \eqref{sing_angles_tet}-\eqref{sing_angles_psi}.
 
The shooting function $ \d{S_{\lambda_1}} $ for the $\lambda_1$-continuation is defined by
\begin{equation*}
\begin{split}
S_{\lambda_1} = \big(
        &  p_{\omega_x }(t_f) , \,\,
            p_{\omega_y} (t_f),\,\,
            p_{\theta} (t_f),\,\,
            p_{\psi} (t_f),\,\,
            p_{\phi} (t_f),\,\,
           H_\gamma(t_f), \\
        &  v_z(t_f) \sin \psi_f + v_y(t_f) \cos \theta_f \cos \psi_f,\,\,
         v_z(t_f) \sin \theta_f - v_x(t_f) \cos \theta_f  \\
         & p_{v_y} \sin \psi_f- (p_{v_x} \sin \theta_f \cos \psi_f+ p_{v_z} \cos \theta_f \cos \psi_f)  \big),
\end{split}
\end{equation*}
where $H_\gamma(t_f)$ with $p^0=-1$ is calculated from \eqref{Hamil_n} and $u_1$ and $u_2$ are given by \eqref{OCPn_u1}.

Note that we can use $S_{\lambda_1}$ as shooting function owing to \OCPReps. For the problem $\MTCP$, if $S_{\lambda_1}=0$, then together with $\omega_x(t_f)=0$ and $\omega_y(t_f)=0$, the final point $(x(t_f),p(t_f))$ of the extremal is then lying on the singular surface and this will cause the failure of the shooting. However, for problem \OCPReps, even when $x(t_f)$ belongs to the singular surface, the shooting problem can still be solved.

Initializing with the solution of $\OCPZ$, we can solve this shooting problem with $\lambda_1=0$, and we get a solution of \OCPReps with the terminal conditions \eqref{inicond_ocpztomtcp}-\eqref{fincond_ocpztomtcp}  (the other states at $t_f$ being free). Then, by continuation, we make $\lambda_1$ vary from $0$ to $1$, and in this way we get the solution of \OCPReps for $\lambda_1=1$. With this solution, we can integrate extremal equations \eqref{sys_atti}, \eqref{sys_orbit_pegasus} and \eqref{sys_adjoint} to get the values of the state variable at $t_f$. We denote $\theta_e := \theta(t_f)$, $\psi_e := \psi(t_f)$, $\phi_e := \phi(t_f)$, $\omega_{xe} := \omega_x(t_f)$ and $\omega_{ye} := \omega_y(t_f)$ the ``natural" conditions obtained at the final time. 

\medskip

In a second step, we use the continuation parameter $\lambda_2$ to act on the final conditions, in order to make them pass from the ``natural" values $\theta_e$, $\psi_e$, $\phi_e$, $\omega_{xe}$ and $\omega_{ye}$, to the desired target values $\theta_f$, $\psi_f$, $\phi_f$, $\omega_{xf}$ and $\omega_{yf}$. The shooting function is
\begin{equation*}
\begin{split}
 S_{\lambda_2} = \big( &
        \omega_x(t_f) - (1-\lambda_2) \omega_{xe} -\lambda_2 \omega_{x_f}, \,
        \omega_y(t_f) - (1-\lambda_2) \omega_{ye} -\lambda_2 \omega_{y_f},\\
       & \theta(t_f)-(1-\lambda_2) \theta_e - \lambda_2 \theta_f ,\,
        \psi(t_f)-(1-\lambda_2) \psi_e - \lambda_2 \psi_f,\,
        \phi(t_f)-(1-\lambda_2) \phi_e - \lambda_2 \phi_f ,\\
       & v_z(t_f) \sin \psi_f + v_y(t_f) \cos \theta_f \cos \psi_f,\,
        v_z(t_f) \sin \theta_f - v_x(t_f) \cos \theta_f ,\, \\
        & p_{v_y} \sin \psi_f- (p_{v_x} \sin \theta_f \cos \psi_f+ p_{v_z} \cos \theta_f \cos \psi_f, \,
        H_\gamma(t_f) \big) .
\end{split}
\end{equation*}
Solving this problem by making vary $\lambda_2$ from $0$ to $1$, we obtain the solution of \OCPReps with the terminal conditions \eqref{OCPc_ic}-\eqref{OCPc_fc}.

\paragraph{$\lambda_3$-continuation}
The parameter $\lambda_3$ is added to the terms induced by the aerodynamic forces in system \eqref{sys_atti}-\eqref{sys_orbit_pegasus}, i.e., replace $c_x$ and $c_z$ by $\lambda_3c_xD$ and $\lambda_3 c_z$. It is the only parameter that acts on the differential system, and it is used to differ the Ariane-type and Pegasus-type rockets.
We let $\lambda_3 = 0$ during the $\lambda_1$-continuation and the $\lambda_2$-continuation, which means during these first two steps of the continuation procedure, we do not consider the influence of the drag and the lift.
Then, if we are dealing with a Pegasus-type rocket, we vary $\lambda_3$ from $0$ to $1$ to derive a solution of problem \OCPReps with drag and lift, and the shooting function remains the same as for the $\lambda_2$-continuation. 
Otherwise, for an Ariane-type rocket, we skip $\lambda_3$-continuation, i.e., we keep $\lambda_3=0$ and go ahead to the $\lambda_4$-continuation. By doing this, the program is suitable for both types of rockets.

\paragraph{$\lambda_4$-continuation}
Finally, in order to compute the solution of $\MTCP$, we use the continuation parameter $\lambda_4$ to pass from \OCPReps to $\MTCP$. We add the parameter $\lambda_4$ to the Hamiltonian $H_{\gamma}$ and to the cost functional \eqref{cost_opcr} as follows:
$$
C_\gamma = t_f + \gamma \int_0^{t_f} (u_1^2 + u_2^2)(1-\lambda_4) \, dt ,
$$
$$
H(t_f,\lambda_4)=\langle p,f(x,u) \rangle +p^0 + p^0\gamma( u_1^2 +u_2^2)(1-\lambda_4).
$$
Then, according to the PMP, the extremal controls are given by $u_i=\mathrm{sat}(-1,u_{ie},1)$, $i=1,2$, where
$$
u_{1e} =   \frac{\bar{b} p_{\omega_y}}{-2 p^0 \gamma (1-\lambda_4) + \bar{b}  \lambda_4 \sqrt{p_{\omega_x}^2+p_{\omega_y}^2} } ,\quad
u_{2e} = \frac{ -\bar{b} p_{\omega_x}}{-2 p^0 \gamma (1-\lambda_4) + \bar{b}  \lambda_4 \sqrt{p_{\omega_x}^2+p_{\omega_y}^2} } .
$$
The shooting function $S_{\lambda_4}$ is the same as $S_{\lambda_2}$, replacing $H_\gamma(t_f)$ with $H_\gamma(t_f,\lambda_4)$. The solution of $\MTCP$ is then obtained by making vary $\lambda_4$ continuously from $0$ to $1$. 

Note again that the above continuation procedure fails in case of chattering, and thus cannot be successful for any possible choice of terminal conditions. In particular, if chattering occurs then the $\lambda_4$-continuation is expected to fail for some value $\lambda_4 = \lambda_4^\ast<1$. But in that case, with this value $\lambda_4^\ast$, we have generated a sub-optimal solution of the problem $\MTCP$, which appears to be acceptable and very interesting in practice. Moreover, the overall procedure is very fast and accurate. Note that the resulting sub-optimal control is continuous.

\medskip

As we have mentioned, it is possible that an optimal trajectory meets the singularities of the Euler angles, i.e., $\cos \psi(t) = 0$, for some $t \in [0,t_f]$. A coordinate system transformation can be made to change the set terminal conditions and avoid meeting the singularity. We will use two numerical tricks as follows. 
The first one is to use a new reference frame $S_R^\ast$, instead of $S_R$, such that the terminal conditions in the new reference frame are easier to solve. 
The new frame $S_R^\ast$ is set by rotating the initial frame $S_R$. This amounts to a nonlinear state transformation, and leads to a preconditioner that makes the proposed continuation procedure more robust.
The second trick is to overcome the singularities of the Euler angles by redefining the vector fields at the singular points.

\subsection{Change of frame}
\label{precond}
We define a new coordinate $S_{R^\prime}$ which is derived by two single-axis rotations from the frame $S_{R}$, i.e., 
$$
S_R  \xrightarrow{R_y(\alpha)} \circ \xrightarrow{R_x(\beta)} S_{R^\prime},
$$
and so the transfer matrix from $S_R$ to $S_R^\ast$ is $L_{R^\prime R} = R_x(\beta)R_z(\alpha)$.
Denoting the Euler angles of $S_b$ with respect to the new frame $S_{R^\prime}$ as $\theta^\prime$, $\psi^\prime$ and $\phi^\prime$, we get the transfer matrix from $S_R^\ast$ to $S_b$ as $L_{bR^\prime} = R_z(\phi^\prime)R_x(\psi^\prime)R_y(\theta^\prime)$.

Therefore, it follows from $L_{bR} L_{R R^\prime }= R_z(\phi)R_x(\psi)R_y(\theta) L_{R^\prime R}^\top= L_{bR^\prime} $ that the angles $\theta^\prime$, $\psi^\prime$ and $\phi^\prime$ are functions of the angles $\theta$, $\psi$, $\phi$, $\alpha$, and $\beta$. 
Moreover, the velocity $\vec{v}$ in the $S_R^\prime$ can also be obtained by $(\vec{v})_R^\prime = L_{R^\prime R}(\vec{v})_R$. The angular velocity $\omega = (\omega_x,\omega_y)$ is expressed in the body frame $S_b$ and it is therefore not altered by this coordinate change.

Given fixed $\alpha$ and $\beta$, the change of frame corresponds to a nonlinear invertible change on the state variable as follows
$$
x^\prime = \mathrm{diag} (L_{R^\prime R} (\vec{v})_R, \varphi_{att}(x), \mathrm{Id}),
$$
where
$\varphi_{att}(\cdot)$ is the mapping from the Euler angles in $S_R$ to the Euler angles in $S_R^\prime$, and $\mathrm{Id}$ is identity matrix of order 3. Note that this transformation is invertible. 
For every given value of $\alpha$, the value of $\beta$ can be chosen such that $\psi^\prime(t_f)=\psi^\prime_f = - \psi^\prime_0=-\psi^\prime(0)$. By doing this, the terminal values on $\psi$ is ``closer" to the origin and hence ``farther" from the singularities of the Euler angles. Indeed, this only reduces the risk of meeting with the singularities of the Euler angles, so we will as well use this numerical trick in the next section to improve the robustness of the numerical procedure.

Further, we explain why this state transformation can be seen as a preconditioner to our numerical method. In view of the shooting method, we use a Newton-like method to solve the boundary value problem. For simplicity of explanation, we use the simplest Newton method and denote the variables of the shooting method as $z$. Given an initial guess of $z=z_0$, the equation $S(z)=0$ is solved iteratively, i.e.,
$$
J(z_{k})z_{k+1} = J(z_{k})z_k - S(z_k).
$$
where $J=\partial{S}/\partial{z}$.
Define a diffeomorphism $\varphi(\cdot)$ such that $z = \varphi(y)$. Then the original problem becomes solving $S(y)=0$ and the Newton method iterative steps become
$$
J(y_{k})y_{k+1} = J(y_{k})y_k - S(y_k),
$$
where 
$$
J(y_{k}) = J(z_{k}) \frac{\partial z}{\partial y}(y_{k}) .
$$
It is easy to see that the matrix $M = (\frac{\partial z}{\partial y} )^{-1}$ here actually acts as a preconditioner in the shooting method and it can be used to get a Jacobian matrix $J$ with a smaller condition number. In the numerical experiments, we use a Fortran subroutine hybrd.f (see \cite{More}) which uses a modification of the Powell hybrid method: the choice of the correction is a convex combination of the Newton and scaled gradient directions, and the updating of the Jacobian by the rank-1 method of Broyden.

Note in addition that the differential equations for the new variable $x^\prime$ keep the same form as the old variable $x$, and by using the PMP, the adjoint vector $p^\prime$ to the new state $x^\prime$ can also be derived from $(x,p)$, i.e., 
$$
p_v^\prime = L_{R^\prime R} p_v, \quad 
(p_{\theta}^\prime,p_{\psi}^\prime,p_{\phi}^\prime)^\top= (\frac{d \varphi_{att}(x)}{dx})^{-\top} (p_{\theta},p_{\psi},p_{\phi})^\top, \quad
(p_{\omega_x}^\prime,p_{\omega_y}^\prime) = (p_{\omega_x},p_{\omega_y}),
$$
and this is also invertible.

To sum up, the reference frame can be chosen such that the problem $\MTCP$ can be easier solved numerically.
However, a priori, we do not know which pair of $(\alpha,\beta)$ is the most suitable to the problem. If we take the pair that makes the condition number of the Jacobian matrix $J$ at the initial time smallest, the value of $|\theta_0^\prime-\theta_f^\prime|$ may become too large (huge pitch angle maneuver) leading to the fail of the numerical strategy and it does not ensure that the condition number remains small in all the four continuation steps. If we choose the pair that makes the value $|\theta_0^\prime-\theta_f^\prime|$ the smallest, it is possible that the Jacobian matrix is ill-conditioned. Therefore, we propose to do the following: 

Frame change heuristics
\begin{itemize}
\label{framechangeheuristics}
\item step1: fix arbitrary $\alpha$, calculate $\beta$ so that $\psi_f=-\psi_0$ and assess the new terminal conditions;
\item step2: solve the problem $\MTCP$ using the continuation procedure proposed;
\item step3: if the step2 succeeds, stop. If the step2 fails, repeat step1-step2 with a new $\alpha=\alpha+\delta \alpha$ ($\delta \alpha$ constant).
\end{itemize}
We will see in the numerical results that this enhances the robustness of our numerical methods though it consumes more time.

\subsection{Singularities of the Euler angles} \label{SingularEuler}
As we can see from \eqref{sys_atti}, when $\psi = \pi/2 + k\pi$, $k\in \mathbb{Z}$, the system encounters singularities due to the fact that the Euler angle $\phi$ is not well defined at these points. Here we call these points singularities of the Euler angles. One way to overcome these singularities is to calculate the limit values of the differential of $\dot{x}$ and $\dot{p}$ at these points. 
Assume that $\dot{\theta}$ is bounded, then we have $(\omega_x \sin \phi + \omega_y \cos \phi) \to 0$ when $\psi \to \pi/2+ k\pi$. Thus we get from
$$
\lim_{\psi \to \pi/2+ k\pi} \dot{\theta} \dot{\phi} = \lim_{\psi \to \pi/2+ k\pi} (\omega_x \sin \phi + \omega_y \cos \phi)^2 \sin \psi =0, \quad
\lim_{\psi \to \pi/2+ k\pi} \dot{\theta}/\dot{\phi} = 1,
$$
that $\dot{\theta} = \dot{\phi} = 0$ when $\psi \to \pi/2+ k\pi$. 
Assume that 
$
\lim_{\psi \to \pi/2+ k\pi}  - \frac{p_{\theta}+p_{\phi} \sin \psi}{\cos \psi} = A < \infty
$,
then it follows from
$$
A = \lim_{\psi \to \pi/2+ k\pi}   - \frac{p_{\theta}+p_{\phi} \sin \psi}{\cos \psi}
=\lim_{\psi \to \pi/2+ k\pi}   \frac{\dot{p}_{\theta}+\dot{p}_{\phi} \sin \psi+p_{\phi} \cos \psi \dot{\psi}}{\sin \psi \dot{\psi}} = -A 
$$
that $A = 0$. Hence, we get
$$
\dot{p}_{\theta} = 0, \quad \dot{p}_{\phi} = 0, \quad \dot{p}_{\psi} = a \sin \theta p_{v x} + a \cos \theta p_{v z},\quad
\dot{p}_{\omega x} = -p_{\psi} \cos \phi,\quad \dot{p}_{\omega y} = p_{\psi} \sin \phi.
$$
Summing up, at points $\psi \to \pi/2+k\pi$, equation \eqref{sys_atti} become
\begin{equation}\label{sys_atti_limit}
\begin{split}
& \dot{\theta}=0 , \quad \dot{\psi}=\omega_x \cos \phi - \omega_y \sin \phi ,\quad \dot{\phi}= 0 ,\quad \dot{\omega }_x= -\bar{b} u_2,\qquad \dot{\omega }_y=  \bar{b} u_1,
\end{split}
\end{equation}
and \eqref{sys_adjoint} become
\begin{equation}\label{sys_adj_limit}
\begin{split}
& \dot{p}_\theta=0 , \quad \dot{p}_\psi= a \sin \theta p_{vx} + a \cos \theta p_{vz},\quad \dot{p}_\phi= 0 ,\quad \dot{p}_{\omega x}= -p_{\psi} \cos \phi,\quad \quad \dot{p}_{\omega y}= p_{\psi} \sin \phi,
\end{split}
\end{equation}
When we are close to such a singularity, we use \eqref{sys_atti_limit} and \eqref{sys_adj_limit} instead of the original system equations \eqref{sys_atti}.

\section{Numerical results} \label{Chp_numerical}

\subsection{Statistical tests}
We test the proposed numerical method in three typical cases of launchers in order to check its robustness with respect to the terminal conditions. 
The first case is the atmospheric ascent phase of the Ariane-type rocket, the second one is the flight phase of the Ariane-type rocket (after the first stage separation), and the third one is the launch phase of the Pegasus-type airborne rocket. The launcher data (see equations \eqref{sys_atti}-\eqref{sys_orbit_pegasus}) in these three cases are given in the Table \ref{table_parameters}. We also indicate the computer features on which the tests were run.
\begin{table} [H]
\centering
\begin{tabular}{|l|c|c|c|}
\hline
& Ariane launch & Ariane flight & Pegasus \\
\hline    
 System parameters   & $a=18$, $\bar{b}=0.0138$ &  $a=25$, $\bar{b}=0.0158$ & $a=24.873$, $\bar{b}=0.0607$ \\
\hline
 Aerodynamical parameters   & $c_x=c_z=0$ &  $c_x=c_z=0$ & $c_x=c_z=5\times10^{-6}$ \\
 \hline
Machine info & \multicolumn{3}{|l|}{CPU: Intel Core i7-4790S CPU 3.2GHz \,\,\,Memory: 7.6 Gio}\\
&  \multicolumn{3}{|l|}{Compiler: gcc version 4.8.4 (Ubuntu 14.04 LTS)} \\
\hline
\end{tabular}
\caption{Launcher data.} 
\label{table_parameters}
\end{table}
We solve the problem $\MTCP$ with different terminal conditions. 
The initial and final conditions are swept in the range given in Table \ref{table_statistic}.
The two last lines of the table define the restrictions applied to the terminal conditions in order to exclude unrealistic cases.
\begin{table} [H]
\centering
\begin{tabular}{|l|c|c|c|}
\hline
	& Ariane launch & Ariane flight & Pegasus \\
\hline
  $v_0$                & $[50,500]\, m/s$  & $[2000,6000]\, m/s$ & $[300,300]\, m/s$\\
  $\theta_{v0}$    & fixed $90^\circ$    & $[0,60]^\circ$      & $[-10,0]^\circ$\\
  $\psi_{v0}$        & fixed $0^\circ$     & fixed $0^\circ$     & fixed $0^\circ$\\
  $\theta_{0}$      & fixed $90^\circ$    & $[0,60]^\circ$      & $[0,50]^\circ$\\  
  $\psi_{0}$          & fixed $0^\circ$     & fixed $0^\circ$     & fixed $0^\circ$\\
  $\theta_{f}$        & $[60,85]^\circ$    & $[0,60]^\circ$      & $[60,90]^\circ$\\
  $\psi_{f}$           & fixed $0^\circ$     & fixed $0^\circ$     & $[0,90]^\circ$\\    
  $\omega_{x 0}$ & fixed $0^\circ$     & fixed $0^\circ$    & $[-10,10]^\circ$\\  
  $\omega_{y 0}$ & fixed $0^\circ$     & fixed $0^\circ$    & fixed $0^\circ$\\   
 \hline    
  $\theta_f-\theta_0$      & free        & $[-20,20]^\circ$           & free \\  
  $\theta_{v0}-\theta_0$ & free        & fixed $0^\circ$            & $[-10,0]^\circ$\\   
\hline
\end{tabular}
\caption{Parameter ranges.} 
\label{table_statistic}
\end{table}
For each variable, we choose a discretization step and solve all the possible combinations of this discretization (factorial experiment). There will be respectively $108$, $234$, $480$ cases for the three launcher configurations.

\paragraph{Results without change of frame}
The statistical results without change of frame are given in Table \ref{table_statistic_result1}.
We see that the success rate for the Ariane launch case (85.4\%) is much higher than the Ariane flight case (50 \%) and the Pegasus case (24.6 \%).

The table presents the average time spent for each continuation stage ($\lambda_1$ to $\lambda_4$) and the average number of simulations performed to achieve each stage. These average performances are assessed on the successful cases.

The main reason for this phenomenon is the singular arc (which is considered as the ``highway" in our problem), which causes chattering. 
This can be seen in the following aspects:
\begin{itemize}
\item The first one is the ratio of parameters $a$ and $b$. 
By comparing the case of Airane launch and the case of Pegasus, we find that when $a/b$ is large (Ariane launch), the continuation procedure works better. 
For a launcher, the parameter $a$ represents the capacity to reorient its velocity, whereas the parameter $b$ represents the capacity to reorient its attitude. 
Hence, when $a/b$ is small (Pegasus), it means that the attitude can be turned quickly to the ``highway", which is the singular surface, and that there will probably be a singular arc in the optimal trajectory. 
\item The second one is the velocity.
In case of Ariane flight, though its ratio $a/b$ is large (even a bit larger than Ariane launch case), the velocity is much bigger than in the Ariane launch case. 
Reorienting the velocity by the same angle requires more time than for a small velocity (see \cite{ZTC1}). That is to say, in order to reorient the velocity vector to the right direction, in case of Ariane flight, it needs to stay more time in the ``highway".
\item The third one is the reorientation angle.
In the case of Pegasus, though the velocity is small, the difference between the initial Euler angles and the final Euler angles are large compared with the two other cases. This means that the velocity also has to have a large change in direction, which leads again the optimal trajectory to have high possibility to contain a singular arc.
\end{itemize}

\begin{table} [H]
\centering
\begin{tabular}{|l|c|c|c|}
\hline    
        \multicolumn{4}{|c|}{Statistical results without change of frame} \\
\hline
  & Ariane launch &Ariane flight  & Pegasus \\        
\hline
  Number of cases                   & $108$ & $234$ & $480$\\  
\hline  
  Rate of success (\%)             &&&\\    
  - Total                                     & $82.4$   & $50.0$   &$24.6$\\   
  - Before $\lambda_4$-continuation & $82.4$ & $75.6$  & $40.2$\\ 
\hline   
  Average execution time (s)    & & &\\                                              
  - Total                                              & $3.44$ & $25.66$  & $26.25$\\   
  - In $\lambda_1$-continuation        & $0.16$ & $2.52$  & $0.60$\\   
  - In $\lambda_2$-continuation         & $2.63$ & $6.30$  & $3.35$\\   
  - In $\lambda_3$-continuation         & $0$      & $0$  & $2.46$\\   
  - In $\lambda_4$-continuation         & $0.62$ & $18.23$  & $14.54$\\  
\hline
  Average number of simulations  & & & \\
  - Total                                                 & $122$ & $195$ & $188$ \\     
  - In $\lambda_1$-continuation           & $20$ & $37$ & $47$ \\   
  - In $\lambda_2$-continuation           & $83$ & $60$ & $37$ \\   
  - In $\lambda_3$-continuation           & $0$   & $0$ & $30$ \\   
  - In $\lambda_4$-continuation           & $19$ & $98$ & $74$ \\                                 
\hline
\end{tabular}
\caption{Statistical result without change of frame.} 
\label{table_statistic_result1}
\end{table}
In term of execution time, it appears that most time is consumed in the $\lambda_4$-continuation (see Ariane flight and Pegasus cases). 
So far, we say that the method is efficient for the Ariane launch case, but it is less efficient and robust for the two other cases.
Now we run the same tests using the change of frame.

\paragraph{Results with change of frame}
The statistical result with change of frame are presented in Table \ref{table_statistic_result2}. 
The change of frame is defined by the angles $(\alpha, \beta)$. The resolution is restarted with varying $(\alpha, \beta)$ values until successful (see \ref{precond} Frame change heuristics). The pair $(\alpha^\ast,\beta^\ast)$ is the pair that makes the continuation procedure works. In Table \ref{table_statistic_result2}, the total execution time includes the $(\alpha, \beta)$ search, while the continuation times account only for the final resolution once $(\alpha^\ast,\beta^\ast)$ have been found.

It is obvious that the proposed continuation procedure is more robust when combined with change of frame. Especially in the Pegasus case, the total success rate is improved from $24.6 \%$ to $65 \%$.

Again, the success rates of Ariane flight and of Pegasus are lower than those of the Ariane launch, and the $\lambda_4$-continuation consumes more time than the three other continuations. 

In fact, when the optimal trajectory does not contain chattering arcs, the total execution time is just about several seconds, however, when the chattering phenomenon tends to appear (by varying $\lambda_4$ from $0$ to $1$), the continuation procedure start to take smaller and smaller steps and the execution time becomes much longer. This is why $\lambda_4$-continuation takes more time. 

Therefore, if we stop the $\lambda_4$-continuation when chattering occurs, we can obtain a sub-optimal solution within a very short time. Moreover, the success rate for obtaining a sub-optimal solution (success rate without $\lambda_4$-continuation) is much higher ($88.9 \%$ for Ariane launch, $88 \%$ for Ariane flight and $83.5 \%$ for Pegasus). This indicates again that most failures happen in the $\lambda_4$-continuation.

We observe that either the solution is obtained very quickly or it can take a long time. Therefore, for a statistical analyses,
we define arbitrarily two group of cases, denoted ``easy" and ``difficult". The easy cases are those solved in less than $50\,s$, whereas the difficult cases are the other ones. 

We see that once the right reference frame is chosen, the time for solving the problem is reasonable even in the difficult cases. Much time is consumed by trying different reference frames. The procedure for searching the best possible values of $\alpha$ and $\beta$, which would ensure convergence, is currently under investigation. The heuristics used for the tests (Frame change heuristics in section \ref{precond}) may certainly be improved taking into account the physical and geometrical features of the problem instance.

\begin{table} [H]
\centering
\begin{tabular}{|l|c|c|c|c|c|c|}
\hline    
        \multicolumn{7}{|c|}{Statistical results with change of frame} \\
\hline
 & \multicolumn {2}{|c|}{Ariane launch} &  \multicolumn {2}{|c|}{Ariane flight}  & \multicolumn {2}{|c|}{Pegasus}  \\        
 \hline 
  Number of cases   & \multicolumn {2}{|c|}{} &  \multicolumn {2}{|c|}{} &  \multicolumn {2}{|c|}{} \\     
  - Total & \multicolumn {2}{|c|}{$108$} & \multicolumn {2}{|c|}{$234$} & \multicolumn {2}{|c|}{$480$}\\ 
  - Easy (total execution time $ < 50\,s$) &  \multicolumn {2}{|c|}{$103$} & \multicolumn {2}{|c|}{$120$} &  \multicolumn {2}{|c|}{$126$}\\
  - Diffi. (total execution time $ \geq 50\,s$) &  \multicolumn {2}{|c|}{$5$} & \multicolumn {2}{|c|}{$114$} &  \multicolumn {2}{|c|}{$354$}\\
\hline  
  Rate of success (\%)    &\multicolumn {2}{|c|}{}&\multicolumn {2}{|c|}{}&\multicolumn {2}{|c|}{}\\  
  - Total & \multicolumn {2}{|c|}{$88.0$}  & \multicolumn {2}{|c|}{$54.7$ }   & \multicolumn {2}{|c|}{$65.0$}\\   
  - Before $\lambda_4$-continuation  & \multicolumn {2}{|c|}{$88.9$}  & \multicolumn {2}{|c|}{$88.0$ }   & \multicolumn {2}{|c|}{$83.5$}\\  
\hline  
  - Average execution time (s)   & Easy & Diffi. & Easy & Diffi. & Easy & Diffi.\\   
   - Total with change of frame    & $2.43$   & $105.70$ &$11.38$ & $930.24$ & $32.6$ & $1122.2$\\                             
   - Total with $(\alpha^\ast,\beta^\ast)$    & $3.73$ & $13.33$ & $10.97$   & $61.45$ & $30.39$ & $50.01$\\   
  \hspace{0.5cm}- In $\lambda_1$-continuation    & $0.22$ & $0.67$  & $0.62$ & $7.78$   & $0.52$ & $2.15$\\   
  \hspace{0.5cm}- In $\lambda_2$-continuation    & $2.63$ & $6.99$  & $4.79$ & $8.00$   & $6.56$ & $6.65$\\   
  \hspace{0.5cm}- In $\lambda_3$-continuation    & $0$    & $0$   & $0$      & $0$      & $4.40$ & $12.23$\\   
  \hspace{0.5cm}- In $\lambda_4$-continuation     & $0.88$ & $5.67$  & $5.56$ & $45.67$ & $18.91$ & $28.98$\\   
\hline
  Average number of simulations  & Easy & Diffi. & Easy  & Diffi. & Easy & Diffi.\\
  - Total with $(\alpha^\ast,\beta^\ast)$        & $123$ & $182$  & $122$& $305$ & $246$ & $387$\\ 
  \hspace{0.5cm}- In the 1st continuation        & $21$  & $25$   & $18$ & $54$  & $21$  & $32$\\   
  \hspace{0.5cm}- In the 2nd continuation        & $83$  & $125$  & $56$ & $55$  & $56$  & $52$\\   
  \hspace{0.5cm}- In the 3rd continuation        & $0$   & $0$    & $0$  & $0$   & $48$  & $150$\\   
 \hspace{0.5cm}- In the 4th continuation         & $20$  & $32$   & $47$ & $196$ & $122$ & $153$\\    
\hline 
  Average times of change of frame & $1.4$ & $12.0$ & $0.0$ & $6.3$  & $0.2$ & $9.7$ \\   
\hline     
\end{tabular}
\caption{Statistic result with change of frame.} 
\label{table_statistic_result2}
\end{table}

Here we test the method with extreme conditions, whereas in real use, only small attitude maneuvers are used for Ariane-type launchers. Among the tested cases of Ariane launch (resp. Ariane flight), there are $86$ cases with maneuver time $t_f < 30\,s$ and the average time to solve the problem is $6.85\,s$ (resp. $3.76\,s$).
For these small maneuvers, the continuation method provides rapidly an optimal control. 

Denote $\delta \psi = |\psi_f-\psi_0|$ the yaw change. We observe from the statistical results that, in case of Pegasus, the average execution time increases when $\delta \psi$ increases. When $\delta \psi \leq 10^\circ$, the average total execution time is $25.07\,s$, and when $\delta \psi \leq 20^\circ$ the total average execution increases to $110.30\,s$. Therefore, when applying to quasi-planar maneuvers ($\delta \psi \leq 10^\circ$) for Pegasus like airborne launchers, the continuation procedure is also efficient.

\subsection{Comparison with the direct approach}
The proposed continuation method is much faster and more precise than the direct approaches. Furthermore, the trajectories obtained, even though they are sub-optimal, are physically feasible because the control remains smooth.
On the other hand, the direct methods may produce oscillatory control if chattering.
In this section, we will use the direct method proposed in \cite{ZTC2}. 
In the direct method, the midpoint rule (2nd order method) is used as discretization method and the number of time steps is set to be 100. Note that in the indirect method, we use Runge-Kutta 4th method which is more precise.

\begin{figure}[H]
\centering
\includegraphics[width=14cm,height=3.5cm]{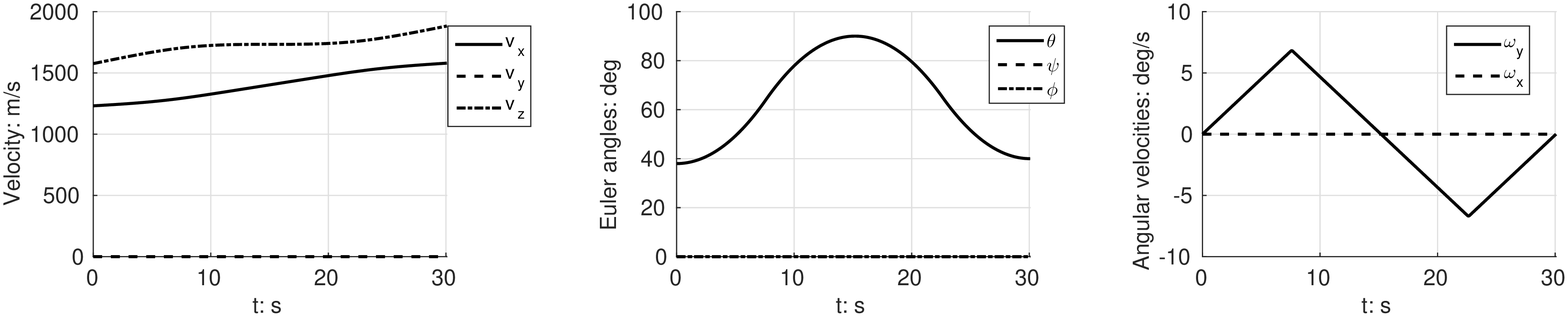}
\caption{State $x(t)$ (continuation method).}
\label{state_easy}
\end{figure}

\begin{figure}[H]
\centering
\includegraphics[width=10cm,height=3.5cm]{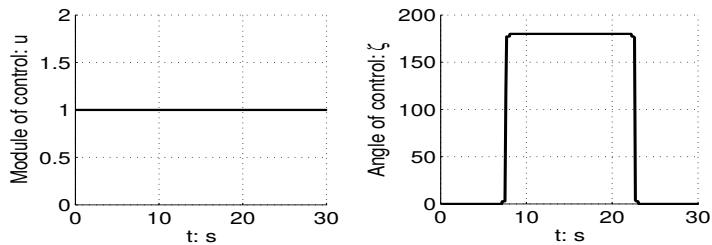}
\caption{Control $u(t)$ (continuation method).}
\label{control_easy}
\end{figure}

\paragraph{Case without chattering}
We report on Figures \ref{state_easy} and \ref{control_easy} an optimal solution of $\MTCP$ obtained from the proposed continuation procedure. 
Here we consider the Ariane flight configuration in a planar case with $\theta_0=38^\circ$, $\theta_f=40^\circ$ and $v_0=2000\,m/s$.
Denote also the module of control as $u$, and thus $u_1 = u \cos \zeta$ and $u_2 = u \sin \zeta$ with $\zeta$ defined as the angle of the control.

For this example, it takes about $8\,s$ to solve the problem using the continuation procedure, while it takes about $40\,s$ using the direct method initialized with an arbitrary constant commands. The optimal solution derived from the direct method is the same with the one reported on the Figures \ref{state_easy} and \ref{control_easy}, except for some small oscillations around zero (namely for the variables $\omega_x$, $\psi$ and $\phi$). Moreover, the maneuver time $t_f=30.03\,s$ given from the direct method is about the same with the maneuver time $t_f=30.07\,s$ obtained from the continuation procedure.

\begin{figure}[H]
\centering
\includegraphics[width=10cm,height=3.5cm]{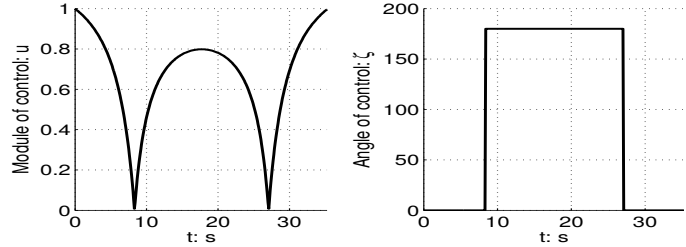}
\caption{Control $u(t)$ (continuation method).}
\label{control_easy_subopti}
\end{figure}

From Figure \ref{control_easy} we see that the obtained optimal control has only two switches (at time $7.5\,s$ and $22.7\,s$ the control angle changes an angle of $\pi$, and thus $u_x$ changes from $+1$ to $-1$ at $7.5\,s$ and changes from $-1$ to $+1$ at $22.7\,s$). 
It is worth noting that for a real launcher, it is better to have a continuous control, since the acceleration of the control is generally limited. For this aim, it suffices to stop the $\lambda_4$-continuation before $1$ to obtain a sub-optimal control. For example, we stop at $0.95$ and get a control in Figure \ref{control_easy_subopti}. This control is feasible on a real launcher and the maneuver time is $t_f=35.25\,s$, five seconds longer than the optimal one.

\medskip

\begin{figure}[H]
\centering
\includegraphics[width=14cm,height=3.5cm]{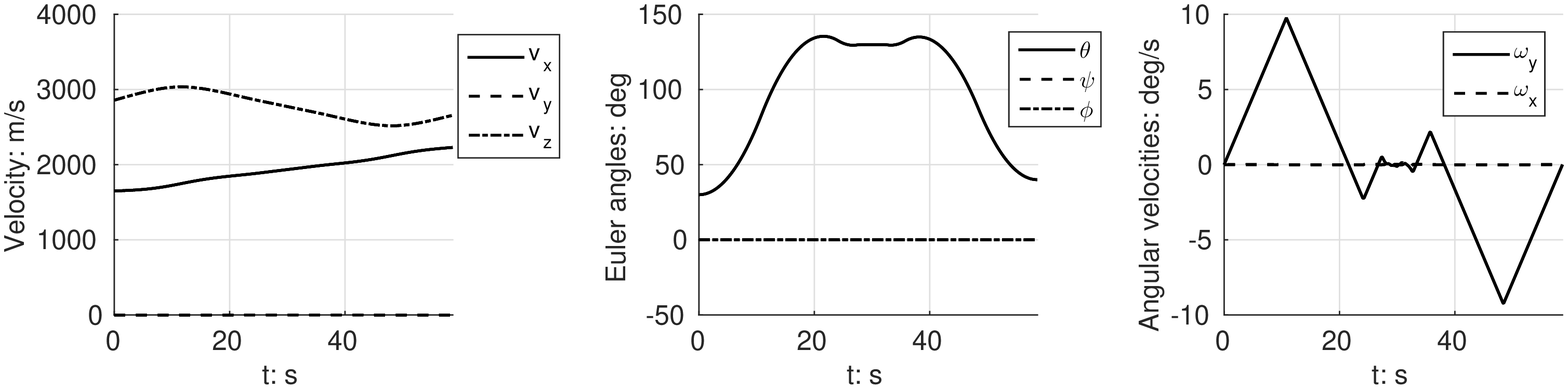}
\caption{State $x(t)$ and control $u(t)$ (direct method).}
\label{state_diffi_bocop}
\end{figure}

\begin{figure}[H]
\centering
\includegraphics[width=10cm,height=3.5cm]{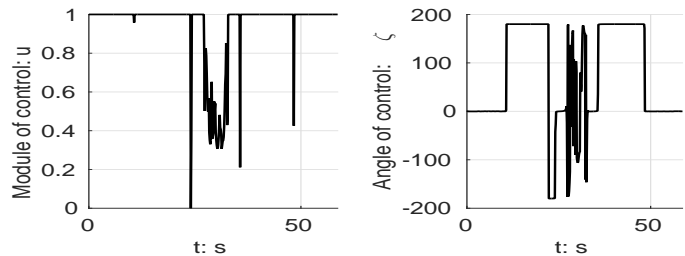}
\caption{Control $u(t)$ (direct method).}
\label{control_diffi_bocop}
\end{figure}

\begin{figure}[H]
\centering
\includegraphics[width=14cm,height=3.5cm]{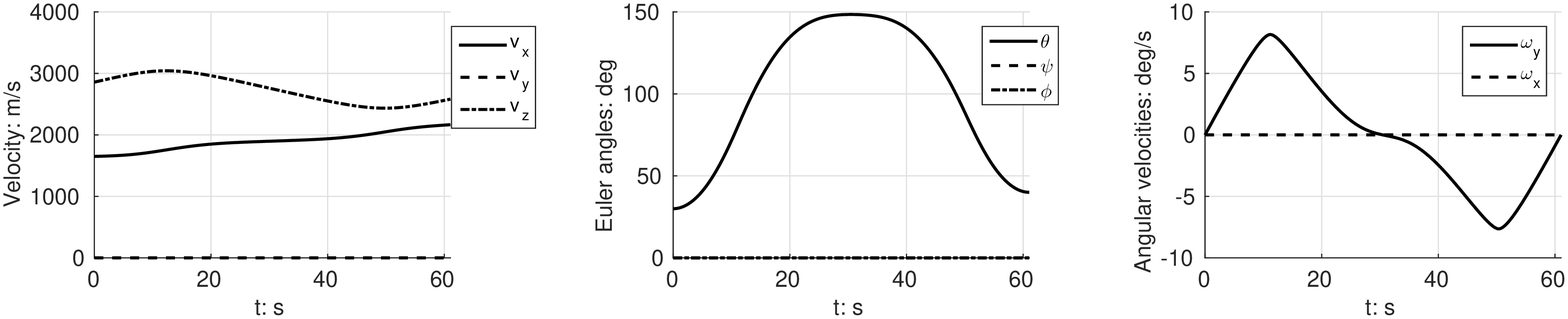}
\caption{State $x(t)$ and control $u(t)$ (continuation method).}
\label{state_diffi}
\end{figure}

\begin{figure}[H]
\centering
\includegraphics[width=10cm,height=3.5cm]{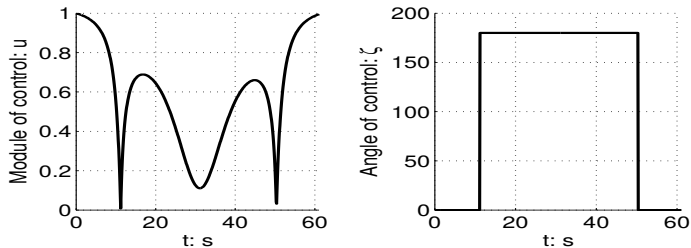}
\caption{Control $u(t)$ (continuation method).}
\label{control_diffi}
\end{figure}

\paragraph{Case with chattering}
In Figures \ref{state_diffi_bocop} and \ref{control_diffi_bocop} a sub-optimal solution obtained from the direct method is illustrated, while in Figures \ref{state_diffi} and \ref{control_diffi} a sub-optimal solution is given by the continuation procedure (parameter $\lambda_4$ stops at $0.994$). Here we still consider the Ariane flight configuration in a planar case with $\theta_0=30^\circ$, $\theta_f=40^\circ$ and $v_0=3200\,m/s$.
When using the continuation method (resp. direct method), the maneuver time obtained is $t_f=61.04\,s$ (resp. $t_f=58.64\,s$) and the problem is solved within $41\,s$ (resp. within $182\,s$). 
We observe that though the maneuver time derived from the direct method is less than with the continuation procedure, the sub-optimal control (see Figure \ref{control_diffi_bocop}) given by the direct method oscillates much with a module less than $1$: this indicates that there is a singular arc in the ``true" optimal trajectory, causing chattering at the junction with regular arcs.
Such solution, with infinitely many control oscillations, cannot be used in the practice. 

From Figure \ref{control_diffi}, we observe that the sub-optimal control obtained from the continuation procedure is continuous and piecewise smooth. In fact, if we observe the time history of the switching function, we will see that the switching function tends to have more switches because of the chattering, and this cause the fail of the continuation method when $\lambda_4>0.995$. 

\medskip

In practice, the attitude maneuvers remain of limited magnitude for a launcher. The chattering phenomenon is thus unlikely to occur, excepted in very specific configurations, for example if the available nozzle deflection is very small.

\section{Conclusion}\label{Chp_conclusion}
In this paper, we have addressed the problem of the minimum time maneuver coupling attitude and trajectory dynamics for a launcher. 
The solution method is based on an indirect approach with a four-stage continuation procedure to solve the optimality conditions derived from the PMP.
Starting from the explicit solution of a simplified problem of order zero, the successive continuations consists of retrieving the true attitude dynamics and terminal conditions. 
A regularization problem is used as an intermediate problem to overcome the problem caused by chattering.
With this procedure, the solution to any problem instance can be obtained from scratch in a quite short time whatever the launcher data and the terminal conditions. 
The minimum time problem is fully solved if possible. In case of a failure due to chattering, a very good sub-optimal solution is available using smooth control law suited to practical application.

The method is exemplified on three representative launcher configurations, respectively for a vertical launch phase, an in-flight phase and an airborne launch phase. 
Statistical tests are run sweeping on the initial and final conditions.
The statistical results show that the developed method is fast and robust. 

The current improvement axes deal with the choice of the coordinate frame that could be chosen from a physical and geometrical considerations. 
Moreover, state constraints on the angular velocity and the attack angle will be considered to provide safer optimal or sub-optimal trajectories for the launchers.

\medskip

\paragraph{Acknowledgment.}
The second author acknowledges the support by FA9550-14-1-0214 of the EOARD-AFOSR.



\begin{thebibliography}{33}
\addcontentsline{toc}{section}{References}

\bibitem{Allgower} E.L. Allgower, K. Georg, 
Numerical continuation methods: an introduction, 
\textit{Springer Science \& Business Media}, 2012.


\bibitem{Betts} J.T. Betts,
Practical methods for optimal control and estimation using nonlinear programming,
Second edition, Advances in Design and Control, 19, 
\textit{Society for Industrial and Applied Mathematics (SIAM)}, Philadelphia, PA, 2010.



 



\bibitem{CHT}
M. Cerf, T. Haberkorn, E. Tr\'elat,
Continuation from a flat to a round Earth model in the coplanar orbit transfer problem,
Optimal Control Appl. Methods, 33 (2012), no. 6, 654--675.

\bibitem{Cesari} L. Cesari,
Optimization - theory and applications. Problems with ordinary differential equations,
Applications of Mathematics, 17, New York: \textit{Springer Verlag}, 1983.


\bibitem{FULLER1} A. T. Fuller, 
An Optimum Non-Linear Control System, 1961,
\textit{In the Proceedings of IFAC Congress}, Moscow, USSR.


\bibitem{Gergaud} J. Gergaud, T. Haberkorn, P. Martinon,
Low thrust minimum fuel orbital transfer: an homotopic approach,
\textit{Journal of Guidance, Control and Dynamics}, 2004, vol. 27, no. 6, p. 1046-1060.










\bibitem{Marchal} C. Marchal, 
Chattering arcs and chattering controls,
\textit{Journal of Optimization Theory and Applications}, 1973, vol. 11, no 5, p. 441-468.

\bibitem{Martinon} P. Martinon, J. Gergaud,
Using switching detection and variational equations for the shooting method,
\textit{Optimal Control Applications and Methods}, 2007, vol. 28, no. 2, p. 95-116.



\bibitem{More} J.J. Mor\'e, D.C. Sorensen, K.E. Hillstrom, et al. 
The MINPACK project. 
\textit{Sources and Development of Mathematical Software}, 1984, p. 88-111.

\bibitem{Pontryagin} L.S. Pontryagin,
Mathematical theory of optimal processes,
\textit{CRC Press}, 1987.




\bibitem{Trelat2} E. Tr\'{e}lat,
Optimal control and applications to aerospace: some results and challenges,
\textit{Journal of Optimization Theory and Applications}, 2012, vol. 154, no. 3, p. 713-758.

\bibitem{TrelatZuazua_JDE2015}
E. Tr\'elat, E. Zuazua,
\textit{The turnpike property in finite-dimensional nonlinear optimal control},
J. Differential Equations {\bf 258} (2015), no. 1, 81--114.

\bibitem{ZELIKIN} M.I. Zelikin, V.F. Borisov, A.J. Krener,
Theory of Chattering Control: with applications to Astronautics, Robotics, Economics, and Engineering,
\textit{Springer}, 1994. 

\bibitem{ZTC1}
J. Zhu, E. Tr\'elat, M. Cerf,
Planar tilting maneuver of a spacecraft: singular arcs in the minimum time problem and chattering,
Preprint arXiv:1504.06219 (2015), 43 pages.

\bibitem{ZTC2}
J. Zhu, E. Tr\'elat, M. Cerf,
Minimum time control of the rocket attitude reorientation associated with orbit dynamics,
Preprint  arXiv:1507.00172 (2015), 33 pages.

\end{thebibliography}
\end{document}